\documentclass{birkmult}

\usepackage[english]{babel}

\usepackage{url}

\begin{document}

\title[Linearization of ODES]{Linearization of ODEs: algorithmic approach}

%----------Author1
\author{V.P.Gerdt}
\address{%
Laboratory of information technologies \\
Joint Institute for Nuclear Research \\
Joliot-Curie 6, 141980 Dubna, Russia
}
\email{gerdt@jinr.ru}

%----------Author2
\author{D.A.Lyakhov}
\address{%
Radiation Gaseous Dynamics Lab \\
Heat and Mass Transfer Institute of the National Academy of Sciences of Belarus\\
P.~Brovka St 15, 220072 Minsk, Belarus
}
\email{lyakhovda@hmti.ac.by}

\date{\empty} %\empty

\maketitle

An ordinary differential equation (ODE) of the first order solved with respect to the derivative it can be always linearized by a suitable point transformation. However, in this case the linearization procedure is not efficient since finding the linearizing transformation is as hard as solving the equation under consideration. For the second order ODEs the situation is different. In this case only equation of the form
\begin{equation}
y''+F_3 (x,y) (y')^3 + F_2 (x,y) (y')^2 + F_1 (x,y) y' + F(x,y)=0.
\label{eq:1}
\end{equation}
may be a candidate for linearization. Sophus Lie \cite{1} designed the following linearizability criterion for equation (\ref{eq:1}): the equation is linerizable by a point transformation if and only if
$$ (F_3)_{xx} - 2 (F_2)_{xy} +(F_1)_{yy} = (3F_1 F_3-F_2^2)_x- 3(FF_3)_y - 3F_3 F_y + F_2 (F_1)_y,$$
$$ F_{yy} - 2(F_1)_{xy} + (F_2)_{xx} = 3(FF_3)_x +(F_1^2-3FF_2)_y + 3F(F_3)_x - F_1 (F_2)_x.$$

A similar criterion was designed for third- and forth-order ODEs \cite{2,3}. By applying the point transformation
$$
u=f(x,y),\quad t=g(x,y),\quad J=f_x g_y - f_y g_x \ne 0
$$
to a linear ODE of the $n$-th order (whose general form is determined by Laguerre's theorem):
$$
u^{(n)} (t) + \sum \limits_{i=0}^{n-3} A_i(t) u^{(i)} (t)=0\,,
$$
we obtain 
\begin{equation}
y^{(n)}(x) + \frac{P(y^{(n-1)},...,y')}{J (g_x+g_y y')^{n-2} }=0\,.
\label{eq:2}
\end{equation}
Here the coefficients of polynomial $P$ are differential polynomials in $f,g$. The formula (\ref{eq:2}) defines the form of an equation to be a candidate for linearization.

Now one asks whether a given ODE of the rational form solved with respect to the highest order derivative
$$
y^{(n)}(x) + \frac{M(y^{(n-1)},...,y' )} {N(y^{(n-1)},...,,y' )} =0
$$
can be linearized by a point transformation. In other words, whether there exist functions $f,g$ such that the equality
$$
\frac{P(y^{(n-1)},...,y')}{J (g_x+g_y y')^{n-2} } = \frac{M(y^{(n-1)},...,y' )} {N(y^{(n-1)},...,,y' )}
$$
of rational functions in $(y^{(n-1)},...,y')$ holds. The last equality is equivalent to the polynomial one
$$
P(y^{(n-1)},...,y') N(y^{(n-1)},...,y' ) - M(y^{(n-1)},...,y' ) J(g_x + g_y y' )^{(n-2)} = 0.
$$
in $(y^{(n-1)},...,y')$.
Since any polynomial is identically zero if and only if its coefficients are zero ones, the linearizability  check is reduced to solvability of the overdetermined system of nonlinear partial differential equations (PDEs). In addition to the last system one more equation must be taken into account. This equation provides dependence of the functions $A_i (t)$ on $t$ only:
$$
(A_i (x,y) )_x g_y - (A_i (x,y) )_y g_x=0.
$$

The consistency of any polynomially-nonlinear PDE system can be verified algorithmically by using the differential Thomas decomposition~\cite{4}. In doing so, the unknown functions are $f,g,A_i$ and their arguments are $(x,y)$.
\vskip 0.5cm
\textbf{Example 1}\,\cite{3}
$$
2 x^2 y y^{(IV)} + x^2 y^2 + 8 x^2 y'y''' + 16 x y y''' + 6 x^2 (y'')^2 + 48 x y'y''
+ 24 y y'' + 24 (y')^2 = 0.
$$

Thomas' decomposition~\cite{4} partition the solution space of this equation into two disjoint  subspaces  defined by the two involutive systems of PDEs. The first subsystem includes 10 equations:
$$
y f_y(x,y) - 2 f(x,y) - 2 f_{xxxx}(x,y) = 0, \quad f_{xxxy} = 0,\quad x^2 f_{xxy} - 2 f_{y} = 0,
$$
$$
x f_{xy} - 2 f_{y} = 0,\quad y f_{yy} - f_{y} = 0, \quad A_0(x,y) (g_x(x,y))^4 - 1 = 0,
$$
$$
g_y(x,y) = 0,\quad (A_0(x,y))_x = 0,\quad (A_0(x,y))_y = 0,\quad A_1(x,y) = 0,
$$
and two inequations:
$$
A_0(x,y) \ne 0, \quad f_y(x,y) \ne 0.
$$

The equations in the system can be readily solved that yields the result:
$$A_0(x,y) = 1,\quad A_1(x,y) = 0,\quad f(x,y) = x^2 y^2,\quad g(x,y)=x.$$
\vskip 0.5cm
\textbf{Example 2\,\cite{3}}. Consider one more ODE of the 4th order
$$
y^{(IV)}(x) + H(x,y) y''(x) + y'(x)^2 = 0.
$$
and verify whether it can be linearized for some kind of  function $H(x,y)$.
For this purpose we include $H(x,y)$ in the list of unknown functions. The Thomas decomposition outputs the empty set. This shows that the equation cannot be linearized for any $H(x,y)$!

The suggested linearizability test is rather simple and efficient. It is fully algorithmic and gives answer  to the question on linearizability of a given ODE by the point transformations.  The PDE system defining the linearizing transformation consists, as a rule, of differentially monomial or binomial equations and can be effectively solved by heuristical methods built-in modern computer algebra systems.

For high-order ODEs, their linearizability is quite exceptional and its verification may need a large volume of symbolic algebraic computation. By this reason it makes sense to throw away apparently inapplicable cases. Thus, to admit linearization, the Lie symmetry algebra of point symmetries for an ODE must have dimension that is not less than the dimension of the symmetry algebra for a linear ODE of the same order. In other words, the dimension of symmetry algebra for the linearizable equation must be strictly higher than the order of the equation. It should be noted that the dimension of Lie symmetry algebra of the  infinitesimal transformations can be algorithmically determined without integration of the determining equations~\cite{5}.
\vskip 0.5cm
\textbf{Example 3}.
$$
y^{(15)}(x) + y'(x)^{15} - 10 y y^{(10)}(x) = 0.
$$
The symmetry Lie algebra of this equation is one dimensional. Therefore this equation cannot be linearized by point transformations. It is remarkable that the second-order ODEs are linearizable if and only if their  symmetry algebra is of maximally possible dimension equal 8~\cite{6}. As to higher order ODEs, their linearizability can be detected by inspection of the abstract Lie symmetry algebra that can also be found without integration of the determining equations~\cite{7}.
\vskip 0.2cm

\textbf{Remark}. Apart from ODEs with polynomial coefficients, the suggested algorithmic approach is also applicable to some cases when the coefficients of an ODE include elementary functions of the independent variable or/and also special functions defined by algebraic differential equations. Furthermore, our approach can be readily generalized to systems of ODEs~\cite{8,9}.

\end{document}